\expandafter\ifx\csname mthreemacsloaded\endcsname\relax\else \fi

\magnification1100
\input amstex


 \catcode`\@=11
 \let\wlog@ld\wlog
 \def\wlog#1{\relax}

 \newif\ifIN@
 \def\m@rker{\m@@rker}
 \def\IN@{\expandafter\INN@\expandafter}
 \long\def\INN@0#1@#2@{\long\def\NI@##1#1##2##3\ENDNI@
    {\ifx\m@rker##2\IN@false\else\IN@true\fi}%
     \expandafter\NI@#2@@#1\m@rker\ENDNI@}
  \newtoks\Initialtoks@  \newtoks\Terminaltoks@
  \def\SPLIT@{\expandafter\SPLITT@\expandafter}
  \def\SPLITT@0#1@#2@{\def\TTILPS@##1#1##2@{%
     \Initialtoks@{##1}\Terminaltoks@{##2}}\expandafter\TTILPS@#2@}
  \newtoks\Trimtoks@

 \def\ForeTrim@{\expandafter\ForeTrim@@\expandafter}
 \def\ForePrim@0 #1@{\Trimtoks@{#1}}
 \def\ForeTrim@@0#1@{\IN@0\m@rker. @\m@rker.#1@%
     \ifIN@\ForePrim@0#1@%
     \else\Trimtoks@\expandafter{#1}\fi}
 
  \def\Trim@0#1@{%
      \ForeTrim@0#1@%
      \IN@0 @\the\Trimtoks@ @%
        \ifIN@
             \SPLIT@0 @\the\Trimtoks@ @\Trimtoks@\Initialtoks@
             \IN@0\the\Terminaltoks@ @ @%
                 \ifIN@
                 \else \Trimtoks@ {FigNameWithSpace}%
                 \fi
        \fi
      }

  \font\titlebold=cmbx12 scaled 1200
  \font\twelvebold=cmbx12
  \font\tenbold=cmbx10
  \font\ninebold=cmbx9
  \font\sevenbold=cmbx7
  \font\fivebold=cmbx5

  \input amssym.def \input amssym
     \font\titlemsa=msam10 at 14.4pt
     \font\titlemsb=msbm10 at 14.4pt
     \font\titleeufm=eufm10 at 14.4pt
     \font\twelvemsa=msam10 scaled 1200
     \font\twelvemsb=msbm10 scaled 1200
     \font\twelveeufm=eufm10 scaled 1200
     \font\ninemsa=msam9
     \font\ninemsb=msbm9
     \font\nineeufm=eufm9

   \ifx\cyrfam\undefined
   \else
     \immediate\write16{}%
     \message{ !!! cyr fonts already defined. !!! }
     \message{ --- edit out superfluous font defs? }
   \fi
   \newfam\cyrfam
       \font\titlecyr=wncyr10 scaled 1440 
       \font\twelvecyr=wncyr10 scaled 1200
       \font\tencyr=wncyr10
       \font\ninecyr=wncyr9
       \font\sevencyr=wncyr7
       \font\sixcyr=wncyr6

   \newfam\eusmfam
       \font\titleeusm=eusm10 scaled 1440
       \font\twelveeusm=eusm10 scaled 1200
       \font\teneusm=eusm10
       \font\nineeusm=eusm9
       \font\seveneusm=eusm7
       
       \font\fiveeusm=eusm5

\let\Cal\cal

    \font\ninemrm=cmr9 
    \font\ninei=cmmi9
    \font\ninesy=cmsy9 
    \skewchar\ninei='177
    \skewchar\ninesy='60

  \font\twelvemrm=cmr10 at 12pt 
  \font\twelvei=cmmi10 at 12pt
  \font\twelvesy=cmsy10 at 12pt

  \font\titlemrm=cmr10 at 14.4pt 
  \font\titlei=cmmi10 at 14.4pt
  \font\titlesy=cmsy10 at 14.4pt


  \def\Smallfonts{\ninepoint}

  \def\Hfont{\titlepoint\bf}
  \def\Authorfont{\twelvepoint\it}
  \def\HHfont{\twelvepoint\bf}
  \def\HHHfont{\bf}
  \def\Bibfont{\tenbf}
  \def\Coordfont{\nineit }

  \def \thfont {\bf }
  \def \pffont {\it\itSpacing }
  \def \rkfont {\bf }
  \def \dffont {\bf }
  \def \egfont {\bf }

 \def\ninepoint{%
  \def\rm{\fam0\ninerm}%
    \textfont0=\ninemrm  \scriptfont0=\sevenrm  \scriptscriptfont0=\fiverm
    \textfont1=\ninei    \scriptfont1=\seveni   \scriptscriptfont1=\fivei
  \def\mit{\fam1\ninei}%
  \def\oldstyle{\fam1\ninei}%
    \textfont2=\ninesy   \scriptfont2=\sevensy  \scriptscriptfont2=\fivesy
    \textfont3=\tenex    \scriptfont3=\tenex    \scriptscriptfont3=\tenex
  \def\it{\fam\itfam\nineit}%
    \textfont\itfam=\nineit
  \def\bf{\ifmmode\fam\bffam\else\ninebf\fi}%
    \textfont\bffam=\ninebold 
    \scriptfont\bffam=\sevenbold 
    \scriptscriptfont\bffam=\fivebold%
  \def\msa{\fam\msafam\ninemsa}%
    \textfont\msafam=\ninemsa 
    \scriptfont\msafam=\sevenmsa
    \scriptscriptfont\msafam=\fivemsa%
  \def\msb{\fam\msbfam\ninemsb}%
    \textfont\msbfam=\ninemsb%
    \scriptfont\msbfam=\sevenmsb%
    \scriptscriptfont\msbfam=\fivemsb%
  \def\eufm{\fam\eufmfam\nineeufm}%
    \textfont\eufmfam=\nineeufm
    \scriptfont\eufmfam=\seveneufm
    \scriptscriptfont\eufmfam=\fiveeufm
   \def\eusm{\fam\eusmfam\nineeusm}%
     \textfont\eusmfam=\nineeusm
     \scriptfont\eusmfam=\seveneusm
     \scriptscriptfont\eusmfam=\fiveeusm
   \def\cyr{\fam\cyrfam\ninecyr}%
     \textfont\cyrfam=\ninecyr
     \scriptfont\cyrfam=\sevencyr
     \scriptscriptfont\cyrfam=\sixcyr
  \setbox\strutbox=\hbox{\vrule
      height7pt depth3pt width0pt}%
   \baselineskip=10.8pt\rm}

 \let\eightpoint\ninepoint 

 \def\tenpoint{%
  \def\rm{\fam0\tenrm}%
    \textfont0=\tenmrm \scriptfont0=\sevenrm \scriptscriptfont0=\fiverm%
  \def\mit{\fam1\teni}%
  \def\oldstyle{\fam1\teni}%
    \textfont1=\teni   \scriptfont1=\seveni  \scriptscriptfont1=\fivei%
    \textfont2=\tensy  \scriptfont2=\sevensy \scriptscriptfont2=\fivesy%
    \textfont3=\tenex  \scriptfont3=\tenex   \scriptscriptfont3=\tenex%
  \def\it{\fam\itfam\tenit}%
    \textfont\itfam=\tenit%
  \def\bf{\ifmmode\fam\bffam\else\tenbf\fi}%
    \textfont\bffam=\tenbold
    \scriptfont\bffam=\sevenbold%
    \scriptscriptfont\bffam=\fivebold%
  \def\msa{\fam\msafam\tenmsa}%
    \textfont\msafam=\tenmsa%
    \scriptfont\msafam=\sevenmsa%
    \scriptscriptfont\msafam=\fivemsa%
  \def\msb{\fam\msbfam\tenmsb}%
    \textfont\msbfam=\tenmsb%
    \scriptfont\msbfam=\sevenmsb%
    \scriptscriptfont\msbfam=\fivemsb%
  \def\eufm{\fam\eufmfam\teneufm}%
   \textfont\eufmfam=\teneufm
   \scriptfont\eufmfam=\seveneufm
   \scriptscriptfont\eufmfam=\fiveeufm
   \def\eusm{\fam\eusmfam\teneusm}%
    \textfont\eusmfam=\teneusm
    \scriptfont\eusmfam=\seveneusm
    \scriptscriptfont\eusmfam=\fiveeusm
   \def\cyr{\fam\cyrfam\tencyr}%
    \textfont\cyrfam=\tencyr
    \scriptfont\cyrfam=\sevencyr
    \scriptscriptfont\cyrfam=\sixcyr
  \setbox\strutbox=\hbox{\vrule %
      height8.5pt depth3.5ptwidth0pt}%
  \baselineskip=\StdBaselineskip\rm}

 \def\twelvepoint{%
  \def\rm{\fam0\twelverm}%
    \textfont0=\twelvemrm \scriptfont0=\tenmrm \scriptscriptfont0=\sevenrm
    \textfont1=\twelvei   \scriptfont1=\teni   \scriptscriptfont1=\seveni
  \def\mit{\fam1\twelvei}%
  \def\oldstyle{\fam1\twelvei}%
    \textfont2=\twelvesy  \scriptfont2=\tensy  \scriptscriptfont2=\sevensy
    \textfont3=\tenex  \scriptfont3=\tenex  \scriptscriptfont3=\tenex
  \def\it{\fam\itfam\twelveit}%
    \textfont\itfam=\twelveit
  \def\bf{\ifmmode\fam\bffam\else\twelvebf\fi}%
    \textfont\bffam=\twelvebold
    \scriptfont\bffam=\tenbold%
    \scriptscriptfont\bffam=\sevenbold%
  \def\msa{\fam\msafam\twelvemsa}%
    \textfont\msafam=\twelvemsa%
    \scriptfont\msafam=\tenmsa%
    \scriptscriptfont\msafam=\sevenmsa%
  \def\msb{\fam\msbfam\twelvemsb}%
    \textfont\msbfam=\twelvemsb%
    \scriptfont\msbfam=\tenmsb%
    \scriptscriptfont\msbfam=\sevenmsb%
  \def\eufm{\fam\eufmfam\twelveeufm}%
   \textfont\eufmfam=\twelveeufm
   \scriptfont\eufmfam=\teneufm
   \scriptscriptfont\eufmfam=\seveneufm
   \def\eusm{\fam\eusmfam\twelveeusm}%
    \textfont\eusmfam=\twelveeusm
    \scriptfont\eusmfam=\teneusm
    \scriptscriptfont\eusmfam=\seveneusm
   \def\cyr{\fam\cyrfam\tencyr}%
    \textfont\cyrfam=\twelvecyr
    \scriptfont\cyrfam=\tencyr
    \scriptscriptfont\cyrfam=\sevencyr
  \setbox\strutbox=\hbox{\vrule
      height10.2pt depth4.55pt width0pt}%
  \baselineskip=14pt\rm}

 \def\titlepoint{%
    \textfont0=\titlemrm \scriptfont0=\twelvemrm \scriptscriptfont0=\tenmrm
    \textfont1=\titlei   \scriptfont1=\twelvei   \scriptscriptfont1=\teni
  \def\mit{\fam1\titlei}%
  \def\oldstyle{\fam1\titlei}%
    \textfont2=\titlesy  \scriptfont2=\twelvesy  \scriptscriptfont2=\tensy
    \textfont3=\tenex
    \scriptfont3=\tenex
    \scriptscriptfont3=\tenex
  \def\it{\fam\itfam\titleit}%
    \textfont\itfam=\titleit
  \def\bf{\ifmmode\fam\bffam\else\titlebf\fi}%
    \textfont\bffam=\titlebold
    \scriptfont\bffam=\twelvebold%
    \scriptscriptfont\bffam=\tenbold%
  \def\msa{\fam\msafam\titlemsa}%
    \textfont\msafam=\titlemsa%
    \scriptfont\msafam=\twelvemsa%
    \scriptscriptfont\msafam=\tenmsa%
  \def\msb{\fam\msbfam\titlemsb}%
    \textfont\msbfam=\titlemsb%
    \scriptfont\msbfam=\twelvemsb%
    \scriptscriptfont\msbfam=\tenmsb%
  \def\eufm{\fam\eufmfam\titleeufm}%
    \textfont\eufmfam=\titleeufm
    \scriptfont\eufmfam=\twelveeufm
    \scriptscriptfont\eufmfam=\teneufm
   \def\eusm{\fam\eusmfam\titleeusm}%
     \textfont\eusmfam=\titleeusm
     \scriptfont\eusmfam=\twelveeusm
     \scriptscriptfont\eusmfam=\teneusm
   \def\cyr{\fam\cyrfam\tencyr}%
    \textfont\cyrfam=\titlecyr
    \scriptfont\cyrfam=\twelvecyr
    \scriptscriptfont\cyrfam=\tencyr
  \setbox\strutbox=\hbox{\vrule
      height12.3pt depth5.54pt width0pt}%
  \baselineskip=16pt\rm}

\newbox\AuthorBox\newbox\TitleBox
\newbox\TFLinebox
\newbox\FLinebox
\newbox\HLinebox
\def\SetTFLinebox#1{\setbox\TFLinebox=\hbox{#1}}
\def\SetFLinebox#1{\setbox\FLinebox=\hbox{#1}}
\def\SetHLinebox#1{\setbox\HLinebox=\hbox{#1}}

 \def\SetAuthorHead#1{%
     \setbox\AuthorBox=\hbox{\ninepoint \it 
           \ignorespaces\frenchspacing#1\unskip}}
 \def\SetTitleHead#1{%
     \setbox\TitleBox=\hbox{\ninepoint \it
           \ignorespaces\frenchspacing#1\unskip}}

  \def\itSpacing{\relax}
  \def\itSpacingOff{\relax}


 \def\Hrule{\hrule width0pt height0pt}

  \newskip\ProcSkip \ProcSkip 8pt plus2pt minus2pt

 \newskip\LastSkip
 \def\SaveLastSkip{\LastSkip\lastskip}
 \def\RestoreLastSkip{\vskip-\LastSkip\vskip\LastSkip}

 \def\NoindentAfter{\everypar={\setbox0=\lastbox\everypar={}}}

 \long\def\H#1\par#2\par{\notenumber=0 \titlepagetrue%
    {
    \baselineskip=20pt
    \parindent=0pt\parskip=0pt\frenchspacing
    \leftskip=0pt plus .2\hsize minus .3\hsize
    \rightskip=0pt plus .2\hsize minus .3\hsize
 \def\\{\unskip\break}%
    \pretolerance=10000 \Hfont #1\unskip\break
     \vskip7pt\Hrule
\hfill \Authorfont #2\hfill\hfill\unskip}
    \vskip48pt plus 4pt minus 4pt
    \par\NoindentAfter\rm}

 \long\def\Hi#1\par#2\par{\notenumber=0 \titlepagetrue%
    {  \baselineskip=0pt  \parindent=0pt\parskip=0pt\frenchspacing
    \leftskip=0pt plus .2\hsize minus .3\hsize
    \rightskip=0pt plus .2\hsize minus .3\hsize
}
    \rm}


 \newdimen\PageRemainder
  \def\SetPageRemainder{
     \PageRemainder=\pagegoal
     \ifdim\PageRemainder=\maxdimen\PageRemainder=\vsize
     \else\advance\PageRemainder by -1\pagetotal\fi}

  \def\Rpt@{}\def\Rpt@@{}

  \long\def\HH#1\par{\par
  \SaveLastSkip\removelastskip\goodbreak
  \ifdim\LastSkip<30pt 
     \LastSkip 30pt
plus 3pt minus 2pt\fi
  \SetPageRemainder\advance\PageRemainder-\LastSkip
  \ifdim\PageRemainder<150pt
       \edef\Rpt@{remain = \the\PageRemainder\noexpand\\
                pagetotal=\the\pagetotal\noexpand\\
                           pagegoal=\the\pagegoal}%
          \fi
   \ifdim\PageRemainder<65pt 
       \ifdim\PageRemainder > 0pt
          \edef\Rpt@@{\noexpand\\
                      Had HH PageRemainder$<$\relax 65pt\noexpand\\
                      Hence forced break!}%
     \vskip 0pt plus .2\PageRemainder\eject 
    \fi\fi
    \vskip\LastSkip\Hrule 
    \pretolerance=10000\rightskip=0pt plus 3em
    \hangafter1 \hangindent=2.2em%
    \noindent
    \HHfont \unskip \Ednote{\Rpt@\Rpt@@}%
            \def\Rpt@{}\def\Rpt@@{}%
            \ignorespaces
            #1\par\rightskip=0pt\pretolerance=\StdPretolerance%
    \NoindentAfter
\tenpoint\rm%
     \medskip \vskip\ProcSkip}

  \long\def\HHH#1\par{\par%
  \SaveLastSkip\removelastskip\goodbreak
  \ifdim\LastSkip<\ProcSkip%
     \LastSkip\ProcSkip\fi
  \SetPageRemainder\advance\PageRemainder-\LastSkip
  \ifdim\PageRemainder<150pt
       \edef\Rpt@{remain = \the\PageRemainder\noexpand\\
                pagetotal=\the\pagetotal\noexpand\\
                           pagegoal=\the\pagegoal}%
       \fi
   \ifdim\PageRemainder<48pt  
        \ifdim\PageRemainder > 0pt
             \edef\Rpt@@{\noexpand\\
                      Had HHH PageRemainder$<$\relax48pt\noexpand\\
                      Hence forced break!}%
       \vskip 0pt plus .2\PageRemainder\eject 
      \fi\fi
   \vskip\LastSkip\par\noindent
   \HHHfont \unskip\Ednote{\Rpt@\Rpt@@}%
  \def\Rpt@{}\def\Rpt@@{}%
  \ignorespaces
   #1\unskip.\quad\rm\ignorespaces
   \ignorepars}

  \long\def\ignorepars#1\par{\def\Test{#1}%
     \ifx\Test\Empty\def\This{\ignorepars}%
        \else\def\This{\Test\par}\fi
           \This}
  \def\Empty{}

 \def\Abstract#1\par{\bgroup\Smallfonts\narrower\HHH #1\par}
 \def\endAbstract{\par\egroup}


 \def\ProcBreak{\par%
    \ifdim\lastskip<8pt%
    \removelastskip%
    \penalty-200\vskip\ProcSkip\fi}

 \def\th#1\par{\ProcBreak \noindent
   {\thfont\ignorespaces
    #1\unskip.}\it\itSpacing\kern.4em\ignorepars}

 \def\endth{\ProcBreak\rm\itSpacingOff }


 \def\pf#1\par{\ProcBreak %
    \noindent\pffont#1\unskip.\rm\itSpacingOff{\kern .7em}\ignorepars}

 \def\endpf{\medskip \ProcBreak } 

  \def\qedbox{\hbox{\vbox{
    \hrule width0.2cm height0.2pt
    \hbox to 0.2cm{\vrule height 0.2cm width 0.2pt
             \hfil\vrule height0.2cm width 0.2pt}
    \hrule width0.2cm height 0.2pt}\kern1pt}}

  \def\qed{\ifmmode\qedbox
    \else\unskip\ \hglue0mm\hfill\qedbox\ProcBreak\fi}

  \def \rk #1\par{\ProcBreak
     \noindent{\rkfont\ignorespaces #1\unskip.}%
     \rm\kern.6em\ignorepars}

  \def \df #1\par{\ProcBreak
     \noindent{\dffont\unskip\ignorespaces #1\unskip.}%
     \rm\kern.6em\ignorepars}

  \def \enddf {\medskip\ProcBreak }

  \def \eg #1\par{\ProcBreak
     \noindent\egfont\unskip\ignorespaces #1\unskip.
     \rm\kern.6em\ignorepars}

  \newdimen\Overhang

   \def\MaxTag@#1#2#3#4#5{\setbox0=\hbox{#4\ignorespaces#2\unskip}%
     \dimen0=\wd0\advance\dimen0 by#3
     \ifdim\dimen0<#5\relax\dimen0=#5\fi
     \expandafter\edef\csname #1Hang\endcsname{\the\dimen0}}

 \def\MaxItemTag#1{\MaxTag@{Item}{#1}{.4em}{\ItemStyle}{\parindent}}%
 \def\MaxItemItemTag#1{%
        \MaxTag@{ItemItem}{#1}{.4em}{\ItemItemStyle}{\parindent}}
 \def\MaxNrTag#1{\MaxTag@{Nr}{#1}{.5em}{\NrStyle}{\parindent}}
 \def\MaxReferenceTag#1{%
        \MaxTag@{Reference}{[#1]}{.6em}{\ninerm}{\parindent}}
 \def\MaxFootTag#1{\MaxTag@{Foot}{#1}{.4em}{\ninerm}{\z@}}

  \def\SetOverhang@{\Overhang=.8\dimen0%
     \advance\Overhang by \wd0\relax
     \ifdim\Overhang>\hangindent\relax
       \advance\Overhang by .25\dimen0%
       \Ednote{Tag is pushing text.}\osumess{Tag is pushing text.}%
     \else\Overhang=\hangindent
     \fi}

   \def\Item#1{\par\noindent
      \hangafter1\hangindent=\ItemHang
      \setbox0=\hbox{\ItemStyle\ignorespaces#1\unskip}%
      \dimen0=.4em\SetOverhang@
      \rlap{\box0}\kern\Overhang\ignorespaces}

   \def\ItemItem#1{\par\noindent
      \hangafter1\hangindent=\ItemItemHang
      \setbox0=\hbox{\ItemItemStyle\ignorespaces#1\unskip}%
      \dimen0=.4em\SetOverhang@
      \advance\hangindent by \ItemHang
      \kern\ItemHang\rlap{\box0}%
      \kern\Overhang\ignorespaces}

  \def\Nr#1{\par\noindent\hangindent=\NrHang 
    \setbox0=\hbox{\NrStyle\ignorespaces#1\unskip}%
    \dimen0=.5em\SetOverhang@
    \rlap{\box0}\kern\Overhang
    \hangindent=\z@\ignorespaces}

   \newskip\Rosterskip\Rosterskip 1pt plus1pt 
   \def\Roster{\par\ifdim\lastskip<\Rosterskip\removelastskip\vskip\Rosterskip\fi
    \bgroup}
   \def\endRoster{\par\global\edef\LastSkip@{\the\lastskip}\removelastskip
       \egroup\penalty-50\LastSkip\LastSkip@\relax
       \ifdim\LastSkip<\Rosterskip\LastSkip\Rosterskip\fi
       \vskip\LastSkip}




 \def\cite#1{
    \def\nextiii@##1,##2\end@{{\frenchspacing\rm 
      \lBr\ignorespaces##1\unskip{\rm,~\ignorespaces##2}\rBr}}%
    \IN@0,@#1@%
    \ifIN@\def\next{\nextiii@#1\end@}\else
    \def\next{{\rm\lBr#1\rBr}}\fi\next}


   \def \Bib#1\par{%
       \par\removelastskip\SetPageRemainder
       \ifdim\PageRemainder < 97pt
        \ifdim\PageRemainder > 0pt
        \vfill\eject
       \fi\fi
    \ProcBreak \par\begingroup\parskip=0 pt%
    \goodbreak \vskip 15 pt plus 10 pt
    \noindent\null\hfill\Bibfont
      \ignorespaces #1\unskip\hfill\null\par 
    \frenchspacing \Smallfonts\rm
    \parskip=2.5 pt plus 1 pt minus.5pt%
    \nobreak\vskip 12pt plus 2pt minus2pt\nobreak
    \leftskip=0 pt \baselineskip=10.5pt}

 \def\ReferenceTagSlide{0em}
  \def\ReferenceTagGap{.5em}

  \def \rf#1{\par\noindent
     \hangafter1\hangindent=\ReferenceHang      
     \setbox0=\hbox{\ninerm[\ignorespaces#1\unskip]}%
     \dimen0=\ReferenceTagGap\SetOverhang@
     \rlap{\kern\ReferenceTagSlide\box0}%
     \kern\Overhang\ignorespaces}

  \def\ref#1\par#2\par#3\par#4\par{%
     \rf{#1}#2\unskip,\ #3\unskip,\
     #4\unskip.}

  \def\endBib{\par\endgroup\vskip 12pt minus 6pt }


  \long\def\Coordinates#1\endCoordinates{
 {\par\vskip4pt\def\\{\unskip, }\Coordfont\baselineskip10.5pt\noindent#1}}

 \def\pagecontents{
  \gdef\Pagetot@l{\pagetotal}
  \ifvoid\TRMargIns\else
    \rlap{\kern\hsize\kern10pt\vbox to 0pt{%
         \box\TRMargIns\vss}}\fi
  \ifvoid\topins\else\unvbox\topins\fi
   \dimen@=\dp\@cclv \unvbox\@cclv 
   \ifvoid\footins\else 
     \vskip\skip\footins
     \footnoterule
     \unvbox\footins\fi
   \ifr@ggedbottom \kern-\dimen@ \vfil \fi}


 \newcount\Ht 

 \def \Acc{\expandafter } 

 \def\swthat{\raise -1.1 ex\hbox{\sam$\widehat{}$}}
 \def\swttilde{\raise -1.2 ex\hbox{\sam$\widetilde{}$}}
 \def \overdot{{\raise .2 ex \hbox to 0pt {\hss\bf\smash{.}\hss}}}
 \def \overcircle{{\raise .1 ex \hbox to 0pt
    {\sam$\eightpoint\scriptstyle\hss\circ\hss$}}}

 \def \Mathaccent#1#2{{\sam 
  \setbox4=\hbox{$\vphantom{#2}$}
  \Ht=\ht4 
  \setbox5=\hbox{${#1}$}
  \setbox6=\hbox{${#2}$}
  \setbox7=\hbox to .5\wd6{}
  \copy7\kern .1\Ht \raise\Ht sp\hbox{\copy5}\kern-.1\Ht
  \copy7\llap{\box6}
  }}

  \def\SwtCheck #1{
        \ifmmode \check{#1}%
                \else \v {#1}%
                \fi}

 \def\barpartial {%
   \kern .17 em
    \overline {\kern -.17 em\partial\kern-.03 em}%
    \kern .03 em}

 
  \def\Overline#1{\setbox1=\hbox{\sam ${#1}$}%
      \ifdim \wd1 > 6pt
    \kern .11 em
    \overline {\kern -.11 em#1\kern-.14 em}
    \kern .14 em
  \else
    \kern .03 em
    \overline {\kern -.03 em#1\kern-.04 em}
    \kern .04 em
  \fi}

 \def\SOverline#1{\setbox1=\hbox{\sam ${#1}$}%
      \ifdim \wd1 > 7pt
    \kern .22 em
    \overline {\kern -.22 em#1\kern-.09 em}%
    \kern .09 em
  \else
    \kern .10 em
    \overline {\kern -.10 em#1\kern-.04 em}%
    \kern .04 em
  \fi}


 \def\Underline#1{\setbox1=\hbox{\sam ${#1}$}%
      \ifdim \wd1 > 6pt
    \kern .11 em
    \underline {\kern -.11 em#1\kern-.14 em}
    \kern .14 em
  \else
    \kern .03 em
    \underline {\kern -.03 em#1\kern-.04 em}
    \kern .04 em
  \fi}

 \def\SUnderline#1{\setbox1=\hbox{\sam ${#1}$}%
      \ifdim \wd1 > 7pt
    \kern .04 em
    \underline {\kern -.04 em#1\kern-.2 em}%
    \kern .2 em
  \else
    \kern .0 em
    \underline {\kern -.0 em#1\kern-.15 em}%
    \kern .15 em
  \fi}


 \def \Blackbox
   {\leavevmode\hskip .3pt \vbox
   {\hrule height 5pt\hbox{\hskip 4.5pt}}\hskip .5pt}

 \def \XX{\Blackbox\kern.5pt\Blackbox} 

  \def\.{.\kern1pt}

    \def\Hyphen{\edef\this{\the\hyphenchar\font}%
          \hyphenchar\font=-1\char\this\hyphenchar\font=\this}

 \ifx\undefined\text
  \def\text#1{\hbox{\rm #1}}\fi 



   \everymath{}  

  \def\PassMath@@{\aftergroup\AfterMath@} 

 \let\PassMath@\PassMath@@

 \def\AfterMath@{\futurelet\next\AfterMathMole@}

 \def\AfterMathMole@{
      \ifcat\next\space
          \def\this{}
      \else
      \ifcat\next\egroup %
        \def\this{\osumess{Handset mathsurround?? ---(see dollar brace)}}%
      \else
      \def\this{\AAfterMath@}
      \fi\fi
      \this}

 \def\hyphen@{-}
 \def\paren@{)}
 \def\apostr@{'}

 \def\MSC#1{\kern-.8\mathsurround#1\kern.8\mathsurround}

 \def\AAfterMath@#1{\def\Next{#1}
    \IN@0\Next @,.;:!?\relax @%
    \ifIN@\def\this{\MSC{\Next}}%
    \else
    \ifx\Next\hyphen@\def\this{\futurelet\next\AfterHyphen@}%
    \else
    \ifx\Next\paren@\def\this{#1}%
    \else 
    \ifx\Next\apostr@\def\this{#1}%
    \else \def\this{\osumess{Handset mathsurround??}%
                 #1}\fi\fi\fi\fi
    \this}

 \def\AfterHyphen@#1{\def\Next{#1}%
   \ifx\Next\hyphen@\def\this{--}\else
   \ifcat\next\space%
   \def\this{\kern-\mathsurround\kern.05em- \Next}\else
   \def\this{\kern-\mathsurround\kern.05em\Hyphen\Next}\fi\fi\this}

 \def\sam{\mathsurround=\z@\let\PassMath@\relax}  %
 \def\mas{\mathsurround=\StdMathsurround\let\PassMath@\PassMath@@}
 
 \def\Mas{\mathsurround=\StdMathsurround
                \everymath{\PassMath@}\let\PassMath@\PassMath@@}

 \def\m@th{\mathsurround=\z@\everymath{}}

 \def\m@@th{\mathsurround=\z@\everymath={}\let\m@th\relax}

\def\underbar#1{$\setbox\z@\hbox{#1}\dp\z@\z@
      \m@th \underline{\box\z@}$\relax}

\def\mathhexbox#1#2#3{\leavevmode
  \hbox{\m@@th$\m@th \mathchar"#1#2#3$}}

\def\dots{\relax\ifmmode\ldots\else$\m@th\ldots\,$\relax\fi}

\def\dotfill{\cleaders\hbox{\m@@th$\m@th \mkern1.5mu.\mkern1.5mu$}\hfill}
\def\rightarrowfill{$\m@th\mathord-\mkern-6mu%
  \cleaders\hbox{\m@@th$\mkern-2mu\mathord-\mkern-2mu$}\hfill
  \mkern-6mu\mathord\rightarrow$\relax}
\def\leftarrowfill{$\m@th\mathord\leftarrow\mkern-6mu%
  \cleaders\hbox{\m@@th$\mkern-2mu\mathord-\mkern-2mu$}\hfill
  \mkern-6mu\mathord-$\relax}

\def\downbracefill{$\m@th\braceld\leaders\vrule\hfill\braceru
  \bracelu\leaders\vrule\hfill\bracerd$\relax}
\def\upbracefill{$\m@th\bracelu\leaders\vrule\hfill\bracerd
  \braceld\leaders\vrule\hfill\braceru$\relax}

\def\angle{{\vbox{\m@@th\ialign{$\m@th\scriptstyle##$\crcr
      \not\mathrel{\mkern14mu}\crcr
      \noalign{\nointerlineskip}
      \mkern2.5mu\leaders\hrule height.34pt\hfill\mkern2.5mu\crcr}}}}

\def\big#1{{\m@@th\hbox{$\left#1\vbox to8.5\p@{}\right.\n@space$}}}
\def\Big#1{{\m@@th\hbox{$\left#1\vbox to11.5\p@{}\right.\n@space$}}}
\def\bigg#1{{\m@@th\hbox{$\left#1\vbox to14.5\p@{}\right.\n@space$}}}
\def\Bigg#1{{\m@@th\hbox{$\left#1\vbox to17.5\p@{}\right.\n@space$}}}
\def\n@space{\nulldelimiterspace\z@ \m@th}

\def\root#1\of{\setbox\rootbox\hbox{\m@@th$\m@th\scriptscriptstyle{#1}$}
  \mathpalette\r@@t}
\def\r@@t#1#2{\setbox\z@\hbox{\m@@th$\m@th#1\sqrt{#2}$\relax}
  \dimen@\ht\z@ \advance\dimen@-\dp\z@
  \mkern5mu\raise.6\dimen@\copy\rootbox \mkern-10mu \box\z@}

\def\mathph@nt#1#2{\setbox\z@\hbox{\m@@th$\m@th#1{#2}$}\finph@nt}

\def\mathsm@sh#1#2{\setbox\z@\hbox{\m@@th$\m@th#1{#2}$}\finsm@sh}

\def\@vereq#1#2{\lower.5\p@\vbox{\m@@th\baselineskip\z@skip\lineskip-.5\p@
    \ialign{$\m@th#1\hfil##\hfil$\crcr#2\crcr=\crcr}}}

\def\mathpalette#1#2{\sam\mathchoice{#1\displaystyle{#2}}%
  {#1\textstyle{#2}}{#1\scriptstyle{#2}}{#1\scriptscriptstyle{#2}}\mas}

\def\widehat#1{\setbox\z@\hbox{\sam$#1$}%
 \ifdim\wd\z@>\tw@ em\mathaccent"0\msbfam@5B{#1}%
 \else\mathaccent"0362{#1}\fi}
\def\widetilde#1{\setbox\z@\hbox{\sam$#1$}%
 \ifdim\wd\z@>\tw@ em\mathaccent"0\msbfam@5D{#1}%
 \else\mathaccent"0365{#1}\fi}

 \def\dots{\relax{}
  \ifmmode\def\thedots{\mdots@}\else\def\thedots{\tdots@}\fi %
  \thedots}

 \let\@ldeqno\eqno\let\@ldleqno\leqno
 \def\eqno{\everymath{}\@ldeqno} \def\leqno{\everymath{}\@ldleqno}

  \let\@ldeqalignno\eqalignno
  \def\eqalignno#1{\sam\@ldeqalignno{#1}\mas}
  \let\@ldeqalign\eqalign
  \def\eqalign#1{\sam\@ldeqalign{#1}\mas}

 \def\overrightarrow#1{\vbox{\m@th\ialign{##\crcr
      \rightarrowfill\crcr\noalign{\kern-\p@\nointerlineskip}
      $\hfil\displaystyle{#1}\hfil$\crcr}}}
 \def\overleftarrow#1{\vbox{\m@th\ialign{##\crcr
      \leftarrowfill\crcr\noalign{\kern-\p@\nointerlineskip}
      $\hfil\displaystyle{#1}\hfil$\crcr}}}
 \def\overbrace#1{\mathop{\vbox{\m@th\ialign{##\crcr\noalign{\kern3\p@}
      \downbracefill\crcr\noalign{\kern3\p@\nointerlineskip}
      $\hfil\displaystyle{#1}\hfil$\crcr}}}\limits}
 \def\underbrace#1{\mathop{\vtop{\m@th\ialign{##\crcr
      $\hfil\displaystyle{#1}\hfil$\crcr\noalign{\kern3\p@\nointerlineskip}
      \upbracefill\crcr\noalign{\kern3\p@}}}}\limits}

  \let\@ldmatrix\matrix
  \let\end@ldmatrix\endmatrix
  \def\matrix{\sam\@ldmatrix}
  \def\endmatrix{\end@ldmatrix\mas}
  \let\@ldgather\gather
  \let\end@ldgather\endgather
  \def\gather{\sam\@ldgather}
  \def\endgather{\end@ldgather\mas}
  \let\@ldalign\align
  \let\end@ldalign\endalign
  \def\align{\sam\@ldalign}
  \def\endalign{\end@ldalign\mas}
  \let\@ldaligned\aligned
  \let\end@ldaligned\endaligned
  \def\aligned{\sam\@ldaligned}
  \def\endaligned{\end@ldaligned\mas}
  \let\@ldtag\tag
  \def\tag{\sam\@ldtag}
   %

   \let\MinCDArrowWidth\minCDaw@




\newskip\insertskipamount\newskip\inserthardskipamount
\insertskipamount 6pt plus2pt 
\inserthardskipamount 6pt
\def\insertskip{\vskip\insertskipamount}
\newcount\SplitTest
\def\SetSplitTest{\SplitTest\insertpenalties
  \insert\topins{\floatingpenalty1}%
  \advance\SplitTest-\insertpenalties}
\def\midinsert{\par
 \SaveLastSkip\penalty-150\SetSplitTest\RestoreLastSkip
 \ifnum\SplitTest=-1
  \@midfalse\p@gefalse\else\@midtrue\fi\@ins}
\def\@ins{\par\begingroup\setbox\z@\vbox\bgroup%
  \vglue\inserthardskipamount}
\def\endinsert{\egroup 
  \if@mid \dimen@\ht\z@ \advance\dimen@\dp\z@
    \advance\dimen@\insertskipamount
    \advance\dimen@\pagetotal\advance\dimen@-\pageshrink
    \ifdim\dimen@>\pagegoal\@midfalse\p@gefalse\fi\fi
  \if@mid%
    \ifdim\lastskip<\insertskipamount\removelastskip\insertskip\fi
    \nointerlineskip\box\z@\penalty-200\insertskip
  \else%
    \SaveLastSkip
    \insert\topins{\penalty100 
    \splittopskip\z@skip
    \splitmaxdepth\maxdimen \floatingpenalty\z@
    \ifp@ge \dimen@\dp\z@
    \vbox to\vsize{\unvbox\z@\kern-\dimen@}
    \else \box\z@\nobreak\insertskip\fi}
    \RestoreLastSkip
   \fi\endgroup}


  \newcount\notenumber
  
  \def\note{\advance\notenumber by 1
    \footnote{\the\notenumber)}}

  \newbox\footbox

  \def\footnote#1{\let\@sf\empty
    \ifhmode\edef\@sf{\spacefactor\the\spacefactor}\/\fi
    \sam${}^{\fam0 #1}$\@sf\vfootnote{#1}}%

  \def\vfootnote#1{\insert\footins\bgroup
     \interlinepenalty100 \splittopskip=1pt
     \floatingpenalty=20000
     \leftskip=0pt\rightskip=0pt%
     \parindent=.3em
     \Smallfonts\rm
     \FootItem@{#1}
     \futurelet\next\fo@t}

  \def\FootItem@#1{\par\hangafter1\hangindent=\FootHang
     \setbox0=\hbox{\ignorespaces#1\unskip}%
     \dimen0=.4em\SetOverhang@
     \noindent\rlap{\box0}\kern\Overhang\ignorespaces}


  \def\fo@t{\ifcat\bgroup\noexpand\next \let\next\f@@t
    \else\let\next\f@t\fi \next}
  \def\f@@t{\bgroup\aftergroup\@foot\let\next}
  \def\f@t#1{\baselineskip=10pt\lineskip=1pt
            \lineskiplimit=0pt #1\@foot}%
  \def\@foot{
        \hbox{\vrule height0pt depth5pt width0pt}
        \egroup}
  \skip\footins=12 pt plus 0pt minus 0pt 
  \count\footins=1000 
  \dimen\footins=8in 



 \def\osumess#1{\EdSpider{\immediate\write16{Line \the\inputlineno: #1}}}%
 \def\HideEdStuff{\gdef\EdSpider##1{}}

 \font\BigSym=cmmi10 scaled \magstep 4

 \def\change{\InLMargin{\hbox{\BigSym \char63\kern10pt}}}

 \def\beginchange{\InLMargin{\hbox{\sam\twelvepoint$\heartsuit$\kern10pt}}}

 \def\endchange{\InLMargin{\hbox{\sam\twelvepoint$\spadesuit$\kern10pt}}}

 \def\InLMargin#1{\strut\vadjust{%
     \kern-\strutdepth
     \vtop to \strutdepth{%
         \baselineskip\strutdepth
         \llap{\sam$\smash{\hbox{\EdSpider{#1}}}$}\null}}}

 \def\strutdepth{\dp\strutbox}
 \def\strutheight{\ht\strutbox}

 \def\NoteInRMargin#1{\strut\vadjust{%
     \kern-1.001\strutdepth
     \vtop to \strutdepth{%
       \baselineskip\strutdepth
       \vss\rlap{\ninepoint\unskip\hskip\hsize
         \vtop to 0pt{%
           \hsize=16em\hfuzz=\hsize
           \leftskip=10pt%
           \rightskip=0pt plus 10000pt%
           \baselineskip=9.8pt\lineskip=.2pt%
           \let\\\break
           \noindent\EdSpider{#1}\vss}%
                \kern10pt}\hbox{}}
       }}

 \def\ednote#1{\NoteInRMargin{\tentt #1}}

 \def\cbar{\InLMargin{%
      \dimen0=\strutdepth\advance\dimen0 by \lineskip
      \vrule width 3pt
      height \strutheight depth \dimen0 \kern
      3pt}}

 \def\ccbar{\InLMargin{%
      \dimen0=2\strutdepth\advance\dimen0 by 2\lineskip
      \vrule width 3pt
        height 3\strutheight depth \dimen0 \kern
      3pt}}

 \newinsert\TRMargIns
 \dimen\TRMargIns=\maxdimen

  \def\Ednote#1{\insert\TRMargIns{%
       \vbox to 0pt{\hsize=140pt\hfuzz=\hsize
           \leftskip=6pt%
           \rightskip=0pt plus 10000pt%
           \baselineskip=9.8pt\lineskip=.2pt%
           \let\\\break
           \SetPageRemainder
           \vglue540pt\vglue-\PageRemainder
           \noindent\EdSpider{\tentt #1}\vss}%
       \smallskip}}

 \def\KillEdStuff{\def\ednote##1{}\def\Ednote##1{}%
      \let\change\relax\let\beginchange\relax\let\endchange\relax
       \let\cbar\relax\let\ccbar\relax}


  \topskip=12pt
  \newskip\StdBaselineskip 
  \StdBaselineskip 12pt
  \lineskip=1.1pt
  \lineskiplimit=.8pt
  \widowpenalty=10000 
  \clubpenalty=10000  
  \abovedisplayskip=6pt plus 1pt minus 1pt
  \abovedisplayshortskip=3pt plus 1.5pt
  \belowdisplayskip=6pt plus 1pt minus 1pt
  \belowdisplayshortskip=5pt plus 1pt minus 1pt
  \hfuzz=1.5pt   

  \def\StdPretolerance{100}
  \tolerance=\StdPretolerance

  \newdimen\StdMathsurround
  \StdMathsurround=1.5pt 
  \mathsurround=\StdMathsurround
  \Mas                   

   \def\prose{\relax\hbox{\kern.6\StdMathsurround}}
  
  \def\StdParskip{0pt}    
  \parskip=\StdParskip
  \parindent=0.5cm
 

  \def\Times{ptmr  } 
  \def\TimesI{ptmri  } 
  \def\TimesB{ptmb  }
  \def\TimesBI{ptmbi  }
  \def\HelveticaN{phvrrn }

  =\Times at 10bp
  =\TimesB at 10bp
  \font\tenit=\TimesI at 10bp
  =\TimesBI at 10bp

  \font\tenmrm=cmr10  


    =\Times at 9bp 
    \font\nineit=\TimesI at 9bp 
    =\TimesB at 9bp 
    =\TimesBI at 9bp 

    =\HelveticaN at 9bp 


  =\Times at 12bp
  \font\twelveit=\TimesI at 12bp
  =\TimesB at 12bp


  \font\titleit=\TimesI at 14.4bp
  =\TimesB at 14.4bp

 \SetAuthorHead{AuthorHead} 
 \SetTitleHead{TitleHead}  


  \def\lBr{\raise.125ex\hbox{[\kern.1125ex}}
  \def\rBr{\raise.125ex\hbox{\kern.1125ex]}}

 \setbox\footbox=\hbox{\Smallfonts 2)~}



  \bgroup
  \catcode`\@=11 
  \gdef\itSpacing{%
     \xspaceskip=.31em plus.1em minus.05em \sfcode `f=2001
     \itWarning@\let\itWarning@\itWarning@@}
  \gdef\itSpacingOff{%
     \xspaceskip=0pt \sfcode `f=1000
     \let\itWarning@\relax}
   \global\let\itWarning@\relax
  \gdef\itWarning@@{\errmessage{%
  Special italic spacing already in force
  (you have probably omitted an ``endth'').
  See itSpacing macro in osuPSfnt.sty
         }}
  \egroup

 \fontdimen1\titlebf=0.0pt
 \fontdimen2\titlebf=3.6135pt
 \fontdimen3\titlebf=2.8908pt
 \fontdimen4\titlebf=1.44539pt
 \fontdimen5\titlebf=6.64882pt
 \fontdimen6\titlebf=14.45398pt
 \fontdimen7\titlebf=1.60439pt

 \fontdimen1\tenbi=0.26794pt
 \fontdimen2\tenbi=2.50937pt
 \fontdimen3\tenbi=2.00749pt
 \fontdimen4\tenbi=1.00374pt
 \fontdimen5\tenbi=4.59717pt
 \fontdimen6\tenbi=10.03749pt
 \fontdimen7\tenbi=1.11415pt

 \fontdimen1\twelverm=0.0pt
 \fontdimen2\twelverm=3.01125pt
 \fontdimen3\twelverm=2.409pt
 \fontdimen4\twelverm=1.2045pt
 \fontdimen5\twelverm=5.39615pt
 \fontdimen6\twelverm=12.045pt
 \fontdimen7\twelverm=1.33699pt

 \fontdimen1\twelveit=0.27731pt
 \fontdimen2\twelveit=3.01125pt
 \fontdimen3\twelveit=2.409pt
 \fontdimen4\twelveit=1.2045pt
 \fontdimen5\twelveit=5.37207pt
 \fontdimen6\twelveit=12.045pt
 \fontdimen7\twelveit=1.33699pt

 \fontdimen1\twelvebf=0.0pt
 \fontdimen2\twelvebf=3.01125pt
 \fontdimen3\twelvebf=2.409pt
 \fontdimen4\twelvebf=1.2045pt
 \fontdimen5\twelvebf=5.5407pt
 \fontdimen6\twelvebf=12.045pt
 \fontdimen7\twelvebf=1.33699pt

 \fontdimen1\tenrm=0.0pt
 \fontdimen2\tenrm=2.50937pt
 \fontdimen3\tenrm=2.00749pt
 \fontdimen4\tenrm=1.00374pt
 \fontdimen5\tenrm=4.49678pt
 \fontdimen6\tenrm=10.03749pt
 \fontdimen7\tenrm=1.11415pt

 \fontdimen1\tenit=0.27731pt
 \fontdimen2\tenit=2.50937pt
 \fontdimen3\tenit=2.00749pt
 \fontdimen4\tenit=1.00374pt
 \fontdimen5\tenit=4.47672pt
 \fontdimen6\tenit=10.03749pt
 \fontdimen7\tenit=1.11415pt

 \fontdimen1\tenbf=0.0pt
 \fontdimen2\tenbf=2.50937pt
 \fontdimen3\tenbf=2.00749pt
 \fontdimen4\tenbf=1.00374pt
 \fontdimen5\tenbf=4.61723pt
 \fontdimen6\tenbf=10.03749pt
 \fontdimen7\tenbf=1.11415pt

 \fontdimen1\ninerm=0.0pt
 \fontdimen2\ninerm=2.25842pt
 \fontdimen3\ninerm=1.80673pt
 \fontdimen4\ninerm=0.90337pt
 \fontdimen5\ninerm=4.0471pt
 \fontdimen6\ninerm=9.03374pt
 \fontdimen7\ninerm=1.00273pt

 \fontdimen1\nineit=0.27731pt
 \fontdimen2\nineit=2.25842pt
 \fontdimen3\nineit=1.80673pt
 \fontdimen4\nineit=0.90337pt
 \fontdimen5\nineit=4.02904pt
 \fontdimen6\nineit=9.03374pt
 \fontdimen7\nineit=1.00273pt

 \fontdimen1\ninebf=0.0pt
 \fontdimen2\ninebf=2.25842pt
 \fontdimen3\ninebf=1.80673pt
 \fontdimen4\ninebf=0.90337pt
 \fontdimen5\ninebf=4.15552pt
 \fontdimen6\ninebf=9.03374pt
 \fontdimen7\ninebf=1.00273pt


 \newcount\MaxSpaceFactor
 \MaxSpaceFactor=3000 

 \def\ItemStyle{\rm}
 \def\NrStyle{\rm}
 \def\ItemItemStyle{\rm}

 \MaxItemTag{(iii)}
 \MaxItemItemTag{(iii)}
 \MaxNrTag{(2)}
 \MaxFootTag{2)}
 \def\ReferenceHang{30pt}

 \catcode`\@=\active


\loadbold

=\Times  
=\Times scaled750
=\Times scaled650
\font\rms=\Times scaled 920 

=\TimesBI scaled 860
=\TimesI scaled 860

\textfont0=\rrm  
\scriptfont0=\erm 
\scriptscriptfont0=\srm

\def\Augment#1#2{%
    \toks0\expandafter{#1}\toks2{#2}%
    \edef#1{\the\toks0\the\toks2}}

 \font\twelverma=\Times  scaled 1200
 \font\tenrma=\Times  scaled 1000
 \font\ninerma=\Times scaled 920
 =\Times scaled 840
 \font\sevenrma=\Times scaled 760
 =\Times scaled 680
 \font\fiverma=\Times scaled 600

 \Augment\tenpoint{%
  \textfont0=\tenrma  \scriptfont0=\sevenrma  
  \scriptscriptfont0=\fiverma  }

 \Augment\ninepoint{%
  \textfont0=\ninerma  \scriptfont0=\sevenrma 
  \scriptscriptfont0=\fiverma}

 \Augment\twelvepoint{%
  \textfont0=\twelverma  \scriptfont0=\ninerma  
  \scriptscriptfont0=\sevenrma}

\mathsurround=1pt
\hsize=13.45truecm
\vsize=19.5truecm
\hoffset=1.25truecm
\voffset=2truecm
\advance\baselineskip by 2pt

\predefine\til{\~}
\def\~#1{\relax\ifmmode\widetilde{#1}\else\til{#1}\fi}

\redefine \le{\leqslant}
\redefine \ge{\geqslant}
\define \wt#1{\mathaccent"0365{#1}}
\define \wh#1{\mathaccent"0362{#1}}

\define \iss{\,\Mathaccent{\raise -.8 ex\hbox{$\widetilde{}$\kern.1em}}\rightarrow\,}

\define \sep{\mathop{\fam0 sep}}

\define \lln{\operatorname{\fam0 log\,}}

\define \coker{\mathop{\fam0 coker}}
\define \kr{\mathop{\fam0 ker}}

\define \Der{\operatorname{\fam0 Der}}

\define \Tr{\operatorname{\fam0 Tr\,}}

\define \Gal{\mathop{\fam0 Gal}}
\define \Hom{\operatorname{\fam0 Hom}}

\Mas
\HideEdStuff
\rm 
 

\def\issn{{\nineit ISSN 1464-8997 (on line) 1464-8989 (printed)}}

\def\gtp{{\nineit Published 10 December 2000: \ \copyright\ Geometry \& 
Topology Publications}}

\def\gtv3{{\nineit Geometry \& Topology Monographs, Volume 3 (2000) --
Invitation to higher local fields}}


\def\lione
{{\rms Geometry \& Topology Monographs}}

\def \litwo{{\rms Volume 3: Invitation to higher local fields
}} 

\def\tinfo #1.#2.#3-#4
{{
\noindent  {\lione} \hfill 
\par 
\vskip-1.5pt
\noindent {\litwo} \hfill
\par 
\vskip-1,5pt
\noindent {\rms Part #1, section #2, pages #3--#4} \hfill
\vskip24pt 
}}

\def\tinfos #1.#2.#3-#4
{{
\noindent  {\lione} \hfill 
\par 
\vskip-1.5pt
\noindent {\litwo} \hfill
\par 
\vskip-1.5pt
\noindent {\rms Pages #3--#4} \hfill
\vskip24pt 
}}

\def\tinfoi #1
{{
\noindent  {\lione} \hfill 
\par 
\vskip-1.5pt
\noindent {\litwo} \hfill
\par 
\vskip-1.5pt
\noindent {\rms Pages iii--xi: Introduction and contents} \hfill
\vskip26pt 
}}


  \def\titlepagehead{\hfil}

  \newif\iftitlepage\titlepagefalse
  \newif\ifblankpage\blankpagefalse
  \def\makeheadline{
     \ifblankpage{}\else%
     \iftitlepage
\vbox{\line{\vbox to 8.5pt{}
\ninerm
\copy\HLinebox \hfill
\hglue5mm\ninebf\folio 
\titlepagehead}}%
      \else
\vbox{\ifodd\pageno\rightheadline\else\leftheadline\fi}%
      \fi\vskip 12pt\fi}%
     \def\rightheadline{\line{\vbox to 8.5pt{}%
      \ninerm
\copy\TitleBox \hfill
\hglue5mm\ninebf\folio}}%
     \def\leftheadline{\line{\vbox to 8.5pt{}%
        \unskip\ninerm\unskip\ninebf\folio\hglue5mm
 \hfill \copy\AuthorBox
}}

 \footline={\ifblankpage{}\else
\iftitlepage\ninepoint\sam\hfill
\line{\vbox to 8.5pt{}
\copy\TFLinebox
\hfill
\hglue5mm 
}
            \else
\ninepoint\sam\hfill
\line{\vbox to 8.5pt{}
\copy\FLinebox
\hfill 
\hglue5mm
}
\hfil\fi\global\titlepagefalse\fi}

\def\blankpage{{\blankpagetrue\noindent\hbox to 10pt{\hss}\vfill
\pagebreak}}

\tenpoint\rm 
 

\pageno=31

\tinfo I.A.31-41

\SetTFLinebox{\gtp }
\SetFLinebox{\gtv3 }
\SetHLinebox{\issn}

\H A. Appendix to Section 2 

Masato Kurihara and Ivan Fesenko

\SetAuthorHead{M. Kurihara and I. Fesenko}
\SetTitleHead{Part I. Appendix to Section 2. \qquad\qquad}

This appendix aims to provide more details on several notions introduced in
section~2, as well as to discuss some basic facts on differentials
and to provide a sketch of the proof of Bloch--Kato--Gabber's theorem.
The work on it was completed after sudden death of Oleg Izhboldin,
the author of section 2.

\HH A1. Definitions and properties of several basic notions 
\unskip\break \phantom{}\enspace  
 (by M. Kurihara)  

Before we proceed to our main topics, we collect here 
the definitions and properties of several basic notions.

\HHH A1.1. Differential modules 

{}
\smallskip

Let $A$ and $B$ be commutative rings such that $B$ is 
an $A$-algebra. 
We define $\Omega_{B/A}^{1}$ to be the $B$-module of regular differentials 
over $A$. 
By definition, this $B$-module $\Omega_{B/A}^{1}$ is a unique 
$B$-module which has the following property. 
For a $B$-module $M$ we denote by $\Der_{A}(B,M)$ the set of 
all $A$-derivations (an $A$-homomorphism $\varphi\colon B \rightarrow 
M$ is called an $A$-derivation if 
$\varphi(xy)=x\varphi(y)+y\varphi(x)$ and 
$\varphi(x)=0$ for any $x \in A$). 
Then, $\varphi$ induces $\Overline{\varphi}\colon \Omega_{B/A}^{1} 
\rightarrow M$ 
($\varphi=\Overline{\varphi} \circ d$ 
where $d$ is the canonical derivation $d\colon B \rightarrow 
\Omega_{B/A}^{1}$), and 
$\varphi \mapsto \Overline{\varphi}$ yields 
an isomorphism 
$$\Der_{A}(B,M) \iss 
\Hom_{B}(\Omega_{B/A}^{1}, M).$$ 

In other words, $\Omega_{B/A}^{1}$ is the $B$-module defined by the
following 

\noindent generators: $dx$ for any $x \in B$

\noindent and relations: 
$$d(xy)=xdy+ydx$$ 
$$dx=0 \quad \text{for any $x \in A$}.$$

\bigskip

If $A={\Bbb Z}$, we simply denote $\Omega_{B/{\Bbb Z}}^{1}$ by 
$\Omega_{B}^{1}$.  

When we consider $\Omega_{A}^{1}$ for a 
local ring $A$, 
the following lemma is very useful. 

\th Lemma  

If $A$ is a local ring, 
we have a surjective homomorphism 
$$A \otimes_{\Bbb Z} A^{*} \longrightarrow \Omega_{A}^{1}$$
$$a \otimes b \mapsto a d\log b=a \frac{db}{b}.$$
The kernel of this map is generated by  elements of the form 
$$\sum_{i=1}^{k}(a_{i} \otimes a_{i}) -\sum_{i=1}^{l}(b_{i} \otimes
b_{i})$$ 
for $a_{i}$, $b_{i} \in A^{*}$ such that 
$\Sigma_{i=1}^{k} a_{i} = \Sigma_{i=1}^{l} b_{i}$. 
\endth

\pf Proof
 
First, we show the surjectivity. 
It is enough to show that $xdy$ is in the image of the above map for 
$x$, $y \in A$. 
If $y$ is in $A^{*}$, 
$xdy$ is the image of $xy \otimes y$. 
If $y$ is not in $A^{*}$, $y$ is in the maximal ideal of $A$, and 
$1+y$ is in $A^{*}$. 
Since $xdy=xd(1+y)$, $xdy$ is the image of $x(1+y) \otimes (1+y)$. 

Let $J$ be the subgroup of $A \otimes A^{*}$ generated by the 
elements 
$$\sum_{i=1}^{k}(a_{i} \otimes a_{i}) -\sum_{i=1}^{l}(b_{i} \otimes
b_{i})$$ 
for $a_{i}$, $b_{i} \in A^{*}$ such that 
$\Sigma_{i=1}^{k} a_{i} = \Sigma_{i=1}^{l} b_{i}$. 
Put $M=(A \otimes_{\Bbb Z} A^{*})/J$. 
Since it is clear that $J$ is in the kernel of the map in the lemma, 
$a \otimes b \mapsto a d \log b$ induces a surjective homomorphism 
$M \rightarrow \Omega_{A}^{1}$, whose injectivity 
we have to show. 

We regard $A \otimes A^{*}$ as an $A$-module via $a(x \otimes y) 
=ax \otimes y$. 
We will show that $J$ is a sub $A$-module of $A \otimes A^{*}$. 
To see this, it is enough to show 
$$\sum_{i=1}^{k}(xa_{i} \otimes a_{i}) -\sum_{i=1}^{l}(xb_{i} \otimes
b_{i}) 
\in J$$ 
for any $x \in A$. 
If $x \not \in A^{*}$, 
$x$ can be written as $x=y+z$ for some $y$, $z \in A^{*}$, 
so we may assume that $x \in A^{*}$. 
Then, 
$$ 
\aligned 
&\sum_{i=1}^{k}(xa_{i} \otimes a_{i}) -\sum_{i=1}^{l}(xb_{i} \otimes b_{i})
\\
&=  \sum_{i=1}^{k}(xa_{i} \otimes xa_{i} - xa_{i} \otimes x)
-\sum_{i=1}^{l}(xb_{i} \otimes xb_{i} - xb_{i} \otimes x) \\
&=  \sum_{i=1}^{k}
(xa_{i} \otimes xa_{i}) -\sum_{i=1}^{l}(xb_{i} \otimes xb_{i}) \in J.
\endaligned 
$$
Thus, $J$ is an $A$-module, and 
$M=(A \otimes A^{*})/J$ is also an $A$-module. 

In order to show the bijectivity of $M \rightarrow \Omega_{A}^{1}$, 
we construct the inverse map 

\noindent$\Omega_{A}^{1} \rightarrow M$. 
By definition of the differential module (see the property after the 
definition), 
it is enough to check that 
the map 
$$
\aligned 
\varphi\colon A \longrightarrow M \qquad & x \mapsto x \otimes x 
\quad 
\text{(if $x \in A^{*}$)} \\ 
&x \mapsto (1+x) \otimes (1+x)
\quad \text{(if $x \not \in A^{*}$)}
\endaligned 
$$
is a ${\Bbb Z}$-derivation. 
So, it is enough to check 
$\varphi(xy)=x\varphi(y)+y\varphi(x)$. 
We will show this  in  the case where both $x$ and $y$ are in the 
maximal ideal of $A$. 
The remaining cases are easier, and are left to the reader. 
By definition, 
$x\varphi(y)+y\varphi(x)$ is the class of 
$$
\aligned 
&x(1+y) \otimes (1+y) + y(1+x) \otimes (1+x) \\
& = (1+x)(1+y) \otimes (1+y) - (1+y) \otimes (1+y)\\ 
&\quad + (1+y)(1+x) \otimes (1+x) - (1+x) \otimes (1+x) \\
&= (1+x)(1+y) \otimes (1+x)(1+y) - (1+x) \otimes (1+x) \\
&\quad - (1+y) \otimes (1+y). 
\endaligned 
$$
But the class of this element in $M$ is the same as the class of 
$(1+xy) \otimes (1+xy)$. 
Thus, $\varphi$ is a derivation. 
This completes the proof of the lemma. 
\qed
\endpf

By this lemma, we can regard $\Omega_{A}^{1}$ as a group defined by the 
following 

\noindent generators: symbols $[a,b\}$ for $a \in A$ and $b \in A^{*}$

\noindent and relations: 
$$\aligned
&[a_{1}+a_{2}, b\}=[a_{1}, b\} +[a_{2}, b\} \\
&[a, b_{1}b_{2}\}=[a, b_{1}\}+[a_{2}, b_{2}\} \\
&\sum_{i=1}^{k}[a_{i}, a_{i}\} = \sum_{i=1}^{l}[b_{i}, b_{i}\}\quad
\text {\rm where $a_{i}$'s and $b_{i}$'s satisfy \quad 
$\sum_{i=1}^{k} a_{i} = \sum_{i=1}^{l} b_{i}$.}
\endaligned 
$$

\HHH A1.2. $n$-th differential forms 

{}
\smallskip

Let $A$ and $B$ be commutative rings such that $B$ is an $A$-algebra. 
For a positive integer $n>0$, 
we define $\Omega_{B/A}^{n}$ by 
$$\Omega_{B/A}^{n}= \bigwedge_{B} \Omega_{B/A}^{1}.$$
Then, $d$ naturally defines an $A$-homomorphism 
$d\colon \Omega_{B/A}^{n} 
\rightarrow \Omega_{B/A}^{n+1}$, and 
we have a complex 
$$ ...\longrightarrow \Omega_{B/A}^{n-1} \longrightarrow \Omega_{B/A}^{n} 
\longrightarrow \Omega_{B/A}^{n+1} \longrightarrow ...$$
which we call the {\it de Rham complex}.

For a commutative ring $A$, which we regard as a ${\Bbb Z}$-module, 
we simply write 
$\Omega_{A}^{n}$ 
for $\Omega_{A/{\Bbb Z}}^{n}$. 
For a local ring $A$, by Lemma A1.1, we have 
$\Omega_{A}^{n}=\bigwedge_{A}^{n}((A \otimes A^{*})/J)$, 
where $J$ is the group as in the proof of Lemma A1.1.
Therefore 
we obtain 

\th Lemma

If $A$ is a local ring, 
we have a surjective homomorphism 
$$
\aligned
A \otimes (A^{*})^{\otimes n} &\longrightarrow \Omega_{A/{\Bbb 
Z}}^{n} \\ 
a \otimes b_{1} \otimes...\otimes b_{n} &\mapsto a \frac{db_{1}}{b_{1}} 
\wedge...\wedge \frac{db_{n}}{b_{n}}.
\endaligned $$
The kernel of this map is generated by  elements of the form 
$$\sum_{i=1}^{k}(a_{i} \otimes a_{i} \otimes b_{1} \otimes...\otimes
b_{n-1}) 
-\sum_{i=1}^{l}(b_{i} \otimes b_{i} \otimes b_{1} \otimes...\otimes
b_{n-1})$$ 
{{\rm(}}where $\Sigma_{i=1}^{k} a_{i} = \Sigma_{i=1}^{l} b_{i}${{\rm)}} 

\noindent and 
$$a \otimes b_{1} \otimes...\otimes b_{n} \quad 
\text{with $b_{i}=b_{j}$ for some $i \neq j$}.$$ 
\endth

\HHH A1.3. Galois cohomology of ${\Bbb Z}/p^{n}(r)$ 
for a field of characteristic $p>0$

{}  
\smallskip

Let $F$ be a field of characteristic $p>0$. 
We denote by $F^{\sep}$ the separable closure of $F$ 
in an algebraic closure of $F$.

We consider Galois cohomology groups $H^{q}(F,-)
:=H^{q}(\Gal(F^{\sep}/F),-)$. 
For an integer $r \ge 0$, 
we define 
$$H^{q}(F, {\Bbb Z}/p(r))=H^{q-r}(\Gal(F^{\sep}/F), 
\Omega_{F^{\sep}, \log}^{r})$$ 
where 
$\Omega_{F^{\sep}, \log}^{r}$ is the logarithmic part of 
$\Omega_{F^{\sep}}^{r}$, namely the subgroup generated by 
$d \log a_{1} \wedge ... \wedge d \log a_{r}$ for all 
$a_{i} \in (F^{\sep})^{*}$.

We have an exact sequence (cf. \cite{I, p.579})
$$ 0 \longrightarrow  
\Omega_{F^{\sep}, \log}^{r} \longrightarrow 
\Omega_{F^{\sep}}^{r} @>{{\bf F}-1}>> 
\Omega_{F^{\sep}}^{r}/d\Omega_{F^{\sep}}^{r-1} \longrightarrow 0$$ 
where ${\bf F}$ is the map 
$${\bf F}(a \frac{db_{1}}{b_{1}} \wedge ...\wedge \frac{db_{r}}{b_{r}}) 
=a^{p} \frac{db_{1}}{b_{1}} \wedge ...\wedge \frac{db_{r}}{b_{r}}.$$ 
Since $\Omega_{F^{\sep}}^{r}$ is an $F$-vector space, 
we have 
$$H^{n}(F, \Omega_{F^{\sep}}^{r})=0$$ 
for any $n >0$ and $r \ge 0$. 
Hence, we also have 
$$H^{n}(F, \Omega_{F^{\sep}}^{r}/d \Omega_{F^{\sep}}^{r-1})=0$$ 
for $n>0$. 
Taking the cohomology of 
the above exact sequence, we obtain 
$$H^{n}(F, \Omega_{F^{\sep}, \log}^{r})=0$$ 
for any $n \ge 2$. 
Further, we have an isomorphism 
$$H^{1}(F, \Omega_{F^{\sep}, \log}^{r})= 
\coker
(\Omega_{F}^{r} @>{{\bf F}-1}>> 
\Omega_{F}^{r}/d\Omega_{F}^{r-1})$$ 
and 
$$H^{0}(F, \Omega_{F^{\sep}, \log}^{r})= \kr
(\Omega_{F}^{r} @>{{\bf F}-1}>>
\Omega_{F}^{r}/d\Omega_{F}^{r-1}).$$

\th Lemma

For a field $F$ of characteristic $p >0$ and $n>0$, we have 
$$H^{n+1}(F, {\Bbb Z}/p\,(n))=
\coker
(\Omega_{F}^{n} @>{{\bf F}-1}>>
\Omega_{F}^{n}/d\Omega_{F}^{n-1})$$ 
and 
$$H^{n}(F, {\Bbb Z}/p\,(n)) = \kr(\Omega_{F}^{n} 
@>{{\bf F}-1}>>
\Omega_{F}^{n}/d\Omega_{F}^{n-1}).$$ 
Furthermore, 
$H^{n}(F, {\Bbb Z}/p\,(n-1))$ is isomorphic to the group which 
has the following 

\noindent generators{{\rm:}} symbols $[a,b_{1},...,b_{n-1}\}$ where 
$a \in F$, and $b_{1}$,...,$b_{n-1} \in F^{*}$ 

\noindent and relations{{\rm:}} 
$$\aligned
&[a_{1}+a_{2}, b_{1},...,b_{n-1}\}=[a_{1}, b_{1},...,b_{n-1}\} 
+[a_{2}, b_{1},...,b_{n-1}\} \\
&[a, b_{1},....,b_{i}b_{i}',...b_{n-1}\}=
[a, b_{1},....,b_{i},...b_{n-1}\} + 
[a, b_{1},....,b_{i}',...b_{n-1}\}\\
&[a, a, b_{2},....,b_{n-1}\} = 0\\
&[a^{p}-a, b_{1}, b_{2},....,b_{n-1}\} = 0\\
&[a, b_{1},....,b_{n-1}\} = 0 \quad \text{where $b_{i}=b_{j}$ 
for some $i \neq j$}.
\endaligned$$
\endth

\pf Proof

 The first half of the lemma 
follows from the computation of 
$H^{n}(F, \Omega_{F^{\sep}, \log}^{r})$ above and 
the definition of 
$H^{q}(F, {\Bbb Z}/p\,(r))$. 
Using 
$$H^{n}(F, {\Bbb Z}/p\,(n-1))=
\coker
(\Omega_{F}^{n-1} @>{{\bf F}-1}>> 
\Omega_{F}^{n-1}/d\Omega_{F}^{n-2})$$
and Lemma A1.2  
we obtain the explicit description of 
$H^{n}(F, {\Bbb Z}/p\,(n-1))$. 
\qed
\endpf

We sometimes use the notation $H_{p}^{n}(F)$ which is defined by 
$$H_{p}^{n}(F) = H^{n}(F, {\Bbb Z}/p\,(n-1)).$$


Moreover, for any $i >1$, we can define ${\Bbb Z}/p^{i}\,(r)$ by using the 
de Rham--Witt complexes instead of the de Rham complex. 
For a positive integer $i>0$, following Illusie \cite{I},  define 
$H^{q}(F, {\Bbb Z}/p^{i}(r))$ by 
$$H^{q}(F, {\Bbb Z}/p^{i}(r)) = H^{q-r}(F, 
W_{i}\Omega_{F^{\sep}, \log}^{r})$$ 
where 
$W_{i}\Omega_{F^{\sep}, \log}^{r}$ is the logarithmic part of 
$W_{i}\Omega_{F^{\sep}}^{r}$.

Though we do not give here the proof, we have 
the following explicit description of 
$H^{n}(F, {\Bbb Z}/p^{i}\,(n-1))$ using the same method as in the case of  
$i=1$. 

\th Lemma

For a field $F$ of characteristic $p >0$ 
let $W_{i}(F)$ denote the ring of Witt vectors of length $i$, and 
let ${\bf F}\colon W_{i}(F) \rightarrow W_{i}(F)$ denote the Frobenius 
endomorphism. 
For any $n>0$ and $i>0$, 
$H^{n}(F, {\Bbb Z}/p^{i}\,(n-1))$ is isomorphic to the group which 
has the following 

\noindent generators{{\rm:}} symbols $[a,b_{1},...,b_{n-1}\}$ where 
$a \in W_{i}(F)$, and $b_{1}$,...,$b_{n-1} \in F^{*}$ 

\noindent and relations{{\rm:}} 
$$\aligned
&[a_{1}+a_{2}, b_{1},...,b_{n-1}\}=[a_{1}, b_{1},...,b_{n-1}\} 
+[a_{2}, b_{1},...,b_{n-1}\} \\
&[a, b_{1},....,b_{j}b_{j}',...b_{n-1}\}=
[a, b_{1},....,b_{j},...b_{n-1}\} + 
[a, b_{1},....,b_{j}',...b_{n-1}\}\\
&[(0,...,0,a,0,...,0), a, b_{2},....,b_{n-1}\} = 0\\
&[{\bf F}(a)-a, b_{1}, b_{2},....,b_{n-1}\} = 0\\
&[a, b_{1},....,b_{n-1}\} = 0 \quad \text{where $b_{j}=b_{k}$ 
for some $j \neq k$}.
\endaligned$$
\endth 

We sometimes use the notation 
$$H_{p^{i}}^{n}(F)= H^{n}(F, {\Bbb Z}/p^{i}\,(n-1)).$$

\HH A2. Bloch--Kato--Gabber's theorem 
\quad  (by I. Fesenko)

{}

For a field $k$ of characteristic $p$ 
denote 
$$
\aligned
&\nu_n=\nu_n(k)=H^{n}(k, {\Bbb Z}/p\,(n))=\kr(\wp\colon \Omega _{k}^{n}
\to \Omega_{k}^{n}/d\Omega _{k}^{n-1}),\\ 
&\wp={\bold F}-1\colon \bigl( a\frac{db_{1}}{b_{1}}\wedge \dots \wedge \frac{db_{n}}{b_{n}}%
 \bigr) \mapsto \left( a^{p}-a\right) \frac{db_{1}}{b_{1}}\wedge \dots \wedge 
\frac{db_{n}}{b_{n}}+ d\Omega _{k}^{n-1}.
\endaligned
$$

Clearly, the image of the differential symbol
$$d_k\colon K_{n}(k)/p \to  \Omega
_{k}^{n}, \qquad 
\left\{ a_{1},\dots ,a_{n}\right\} 
\mapsto \frac{da_{1}}{a_{1}}\wedge
\dots \wedge \frac{da_{n}}{a_{n}}
 $$ 
is inside $\nu_n(k)$.
We shall sketch the proof of  Bloch--Kato--Gabber's theorem
which states that $d_k$ is an isomorphism between
$K_{n}(k)/p$ and $\nu_n(k)$.

\HHH A2.1. Surjectivity of the differential symbol
$d_k\colon K_{n}(k)/p \to  \nu_n(k)$

\phantom{}\smallskip \par 

It seems impossible to suggest a shorter proof
than  original Kato's proof in \cite{K, \S1}. 

We can argue by induction on $n$; the case of  $n=1$ is obvious,
so assume $n>1$.

\df Definitions--Properties 

\phantom{}\par
\Roster

\Item{(1)} Let $\{b_i\}_{i\in I}$ be a $p$-base of $k$
($I$ is an ordered set).
Let $S$ be the set of all  strictly increasing maps 
$$s\colon \{1,\dots,n\}\to I.$$
For two maps $s,t\colon \{1,\dots,n\}\to I$
write $s<t$ if $s(i)\le t(i)$ for all $i$
and $s(i)\not= t(i)$ for some $i$.

\Item{(2)} Denote $d\lln a:= a^{-1}da$.
Put  
$$\omega_s=d\lln b_{s(1)}\wedge\cdots\wedge d\lln b_{s(n)}.$$
Then $\{\omega_s:s\in S\}$ is a basis of $\Omega_k^n$ over $k$.

\Item{(3)} For a map $\theta\colon I\to \{0,1,\dots, p-1\}$
such that $\theta(i)=0$ for almost all $i$  set 
$$b_{\theta}=\prod
b_i^{\theta(i)}.$$
Then $\{b_{\theta}\omega_s\}$ is a basis of $\Omega_k^n$ over $k^p$. 

\Item{(4)}Denote by $\Omega_k^n(\theta)$ the $k^p$-vector space generated by
$b_{\theta}\omega_s, s\in S$.
Then $\Omega_k^n(0)\cap d\Omega_k^{n-1}=0$.
For  an extension $l$ of $k$, such that $k\supset l^p$, 
denote by $\Omega_{l/k}^n$ 
the module of 
relative differentials.
Let $\{b_i\}_{i\in I}$ be a $p$-base of $l$ over $k$.
Define $\Omega_{l/k}^n(\theta)$ for a map $\theta\colon I\to \{0,1,\dots, p-1\}$
 similarly to the previous definition. 
The cohomology group of the complex 
$$\Omega_{l/k}^{n-1}(\theta)\to
\Omega_{l/k}^n(\theta)\to \Omega_{l/k}^{n+1}(\theta)$$
is zero if $\theta\not=0$ and is $\Omega_{l/k}^n(0)$ if $\theta=0$.
\endRoster 
\enddf

We shall use {\it Cartier's theorem}
(which can be more or less easily proved by induction on $|l:k|$): 
the sequence
$$0\to l^*/k^*\to \Omega_{l/k}^1\to \Omega_{l/k}^1/dl$$
is exact, where the second map
is defined as $b\mod k^*\to d\lln b$ and the third map
is the map 
$ad\lln b\mapsto (a^p-a)d\lln b + dl$.

\th Proposition

Let $\Omega_k^n(<\!\! s)$ be the $k$-subspace of $\Omega_k^n$ 
generated by all
$\omega_t$ for $s>t\in S$.

Let $k^{p-1}=k$ and let $a$ be a non-zero element of $k$.
Let $I$ be finite.
Suppose that $$(a^p-a)\omega_s\in \Omega_k^n(<\!\! s)+d\Omega_k^{n-1}.$$

Then there are $v\in \Omega_k^n(<\!\! s)$
and $$x_i\in k^p(\{b_j:j\le s(i)\})\quad\text{\it for}\quad 1\le i\le n$$ such that
$$a\omega_s=v+d\lln x_1\wedge\cdots\wedge d\lln x_n.$$
\endth

\pf Proof of the surjectivity of the differential symbol

First, suppose that $k^{p-1}=k$ and 
 $I$ is finite. 
Let $S=\{s_1,\dots, s_m\}$ with $s_1>\cdots>s_m$.
Let $s_0\colon \{1,\dots, n\}\to I$ be a map 
such that $s_0>s_1$.
Denote by $A$ the subgroup of $\Omega_k^n$ generated by
$d\lln x_1\wedge\cdots\wedge d\lln x_n$.
Then $A\subset \nu_n$. By induction on $0\le j\le m$ using the proposition
it is straightforward to show that $\nu_n\subset A+ \Omega_k^n(<\!\! s_j)$, and hence
$\nu_n=A$.

To treat the general case 
put $c(k)=\coker(k_n(k)\to \nu_n(k))$.
Since every field is the direct limit of finitely generated fields
and the functor $c$ commutes with direct limits,
it is sufficient to show that $c(k)=0$ for
a finitely generated field $k$.
In particular, we may assume that $k$ has a finite $p$-base.
For a finite extension  $k'$ of $k$ there is a commutative diagram 
$$
\CD
k_{n}(k') @>>>  \nu_n({k'})\\
@VN_{k'/k}VV @V \Tr_{k'/k}VV \\
k_{n}(k) @>>>  \nu_n(k).
\endCD
$$
Hence the composite $c(k)\to c(k')@>\Tr_{k'/k}>> c(k)$
is multiplication by $|k':k|$.
Therefore, if $|k':k|$ is prime to $p$ then $c(k)\to c(k')$ is injective.

Now pass from $k$ to a field $l$ which is the compositum of all $l_i$
where $l_{i+1}=l_i(\root{p-1}\of{l_{i-1}})$, $l_0=k$.
Then $l=l^{p-1}$.
Since $l/k$ is separable, $l$ has a finite $p$-base
and by the first paragraph of this proof
$c(l)=0$.
The degree of every finite subextension in $l/k$ is prime to $p$,
and by the second paragraph of this proof we conclude $c(k)=0$, as required.
\qed
\endpf

\pf Proof of Proposition

First we prove the following lemma which will help us later for 
fields satisfying $k^{p-1}=k$ to choose a specific $p$-base of $k$. 

\th Lemma

Let $l$ be a purely inseparable extension of $k$
of degree $p$ and let $k^{p-1}=k$. Let
$f\colon l\to k$ be a $k$-linear map. Then 
there is a non-zero $c\in l$ such that
$f(c^i)=0$ for all $1\le i\le p-1$.
\endth

\pf Proof of Lemma

The $l$-space of $k$-linear maps from $l$ to $k$
is one-dimensional, hence $f=ag$ for some $a\in l$, where
$g\colon l=k(b)\to \Omega_{l/k}^1/dl\iss k$,
 $x\mapsto xd\lln b$ mod  $dl$ for every $x\in l$.
Let $\alpha= gd\lln b$ generate the one-dimensional space
$\Omega_{l/k}^1/dl$ over $k$.
Then there is $h\in k$ such that $g^pd\lln b-h\alpha\in dl$.
Let $z\in k$ be such that $z^{p-1}=h$. Then $((g/z)^p-g/z)d\lln b\in dl$ and
by  
Cartier's theorem  we deduce that there is $w\in l$ such that
$(g/z)d\lln b=d\lln w$. Hence $\alpha=zd\lln w$ and
$\Omega_{l/k}^1=dl\cup kd\lln l$.

If $f(1)=ad\lln b\not=0$, then $f(1)=gd\lln c$ with $g\in k, c\in l^*$
and hence $f(c^i)=0$ for all $1\le i\le p-1$.
 \qed
\endpf

Now 
for $s\colon\{1,\dots,n\}\to I$ as in the statement of the Proposition
denote 
$$k_0=k^p(\{b_i: i<\!\! s(1)\}), \quad k_1=k^p(\{b_i: i\le s(1)\}),
\quad 
k_2=k^p(\{b_i: i\le s(n)\}).$$
 Let $|k_2:k_1|= p^r$.

Let $a=\sum_\theta x_{\theta}^p b_\theta$.
Assume that $a\not\in k_2$. Then let $\theta, j$ be such that
$j>s(n)$ is the maximal index for which
$\theta(j)\not=0$ and $x_{\theta}\not=0$. 

$\Omega_k^n(\theta)$-projection of $(a^p-a)\omega_s$  is
equal to  
$-x_{\theta}^p
b_{\theta}\omega_s\in \Omega_k^n(<\!\! s)(\theta)+d\Omega_k^{n-1}(\theta)$.
Log differentiating, we get
$$-x_{\theta}^p
\bigl(\sum_i \theta(i)d\lln b_i\bigr)b_\theta\wedge\omega_s\in d\Omega_k^n(<\!\! s)(\theta)$$
which contradicts 
$-x_{\theta}^p \theta(j)b_\theta d\lln b_j\wedge \omega_s\not\in d\Omega_k^n(<\!\! s)(\theta)$. Thus, $a\in k_2$.

Let $m(1)<\dots< m(r-n)$ be integers such that
the union of $m$'s and $s$'s is equal to $[s(1),s(n)]\cap \Bbb Z$.
Apply the Lemma to the linear map
$$
f\colon k_1\to \Omega_{k_2/k_0}^r/d\Omega_{k_2/k_0}^{r-1}\iss k_0,
\quad b\mapsto ba\omega_s\wedge d\lln b_{m(1)}\wedge\cdots\wedge
d\lln b_{m(r-n)}.
$$
Then there is a non-zero $c\in k_1$ such that 
$$
c^ia\omega_s\wedge d\lln b_{m(1)}\wedge\cdots\wedge
d\lln b_{m(r-n)}\in d\Omega_{k_2/k_0}^{r-1}\quad\text{ for $1\le i\le p-1$}.
$$
Hence  $\Omega_{k_2/k_0}^{r}(0)$-projection of $c^ia\omega_s\wedge d\lln b_{m(1)}\wedge\cdots\wedge d\lln b_{m(r-n)}$ for $1\le i\le p-1$ is zero.

If $c\in k_0$ then $\Omega_{k_2/k_0}^{r}(0)$-projection of $a\omega_s\wedge d\lln b_{m(1)}\wedge\cdots\wedge d\lln b_{m(r-n)}$ is zero.
Due to the definition of $k_0$ we get 
$$\beta=(a^p-a)\omega_s\wedge d\lln b_{m(1)}\wedge\cdots\wedge d\lln b_{m(r-n)}\in d\Omega_{k_2/k_0}^{r-1}.$$
Then  $\Omega_{k_2/k_0}^{r}(0)$-projection of $\beta$ is zero, and so is
$\Omega_{k_2/k_0}^{r}(0)$-projection of $$a^p\omega_s\wedge d\lln b_{m(1)}\wedge\cdots\wedge d\lln b_{m(r-n)},$$ a contradiction.
Thus, $c\not\in k_0$.

From $dk_0\subset \sum _{i<s(1)} k^pdb_i$ we deduce $dk_0\wedge \Omega_k^{n-1}\subset
\Omega_k^n(<\!\! s)$.
Since $k_0(c)=k_0(b_{s(1)})$, there are $a_i\in k_0$ such that
$b_{s(1)}=\sum_{i=0}^{p-1} a_ic^i$.
Then 
$$ad\lln b_{s(1)}\wedge\cdots\wedge d\lln b_{s(n)}
\equiv a' d\lln b_{s(2)}\cdots\wedge d\lln b_{s(n)}\wedge
d\lln c \mod \Omega_k^n(<\!\! s).$$

Define $s'\colon \{1,\dots,n-1\}\to I$ by $s'(j)=s(j+1)$.
Then 
$$a\omega_s=v_1+a'\omega_{s'}\wedge d\lln c\quad
\text{ with $v_1\in \Omega_k^n(<\!\! s)$}
$$
 and
$c^ia'\omega_{s'}\wedge d\lln c\wedge  d\lln b_{m(1)}\wedge\cdots\wedge d\lln b_{m(r-n)}\in d\Omega_{k_2/k_0}^{r-1}$.
The set $$I'=\{c\}\cup\{b_i:s(1)<i\le s(n)\}$$ is a $p$-base of $k_2/k_0$.
Since $c^ia'$ for $1\le i\le p-1$ have zero $k_2(0)$-projection
with respect to $I'$, there are $a_0'\in k_0$, $a_1'\in \oplus_{\theta\not=0} k_1b'_{\theta}$ with $b'_{\theta}=\prod_{s(1)<i\le s(n)}b_i^{\theta(i)}$
such that $a'=a_0'+a_1'$. 

The image of $a\omega_s\wedge
  d\lln b_{m(1)}\wedge\cdots\wedge d\lln b_{m(r-n)}$ with respect to the Artin--Schreier map belongs to
$\Omega_{k_2/k_0}^r$ and so is
$$({a'}^p-a')d\lln c\wedge \omega_{s'}\wedge  d\lln b_{m(1)}\wedge\cdots\wedge d\lln b_{m(r-n)}$$ which is the image of $$a'd\lln c\wedge \omega_{s'}\wedge  d\lln b_{m(1)}\wedge\cdots\wedge d\lln b_{m(r-n)}.$$
Then ${a'}^p-a_0'$, as $k_0(0)$-projection of ${a'}^p-a'$, is zero.
So $a'-{a'}^p=a_1'$.

Note that $d(a_1'\omega_{s'})\wedge d\lln c\in d\Omega_{k/k_0}^n(<\!\! s)=d\Omega_{k/k_0}^{n-1}(<\!\! s) \wedge
d\lln c$. 

\noindent Hence $d(a_1'\omega_{s'})\in d \Omega_{k/k_0}^{n-1}(<\!\! s)+
d\lln c\wedge d\Omega_{k/k_0}^{n-2}$. Therefore  
$d(a_1'\omega_{s'})\in d \Omega_{k/k_1}^{n-1}(<\!\! s)$ and
$a_1'\omega_{s'}=\alpha+\beta$ with $\alpha\in  \Omega_{k/k_1}^{n-1}(<\!\! s)$,
$\beta\in \kr (d\colon \Omega_{k/k_1}^{n-1}\to \Omega_{k/k_1}^{n})$.

Since $k(0)$-projection of $a_1'$ is zero, $\Omega_{k/k_1}^{n-1}(0)$-projection of
$a_1'\omega_{s'}$ is zero. Then we deduce that 
$\beta(0)=\sum_{x_t\in k_1, t<s'} x_t^p\omega_t$, so
$a_1'\omega_{s'}=\alpha+\beta(0)+(\beta-\beta(0))$.
Then $\beta-\beta(0)\in
\kr (d\colon \Omega_{k/k_1}^{n-1}\to \Omega_{k/k_1}^{n})$, 
so 
$\beta-\beta(0)\in d\Omega_{k/k_1}^{n-2}$.
Hence   $(a'-{a'}^p)\omega_{s'}=a_1'\omega_{s'}$ belongs to $\Omega_{k/k_1}^{n-1}(<\!\! s')+d\Omega_k^{n-2}$.
By induction on $n$, there are $v'\in \Omega_k^{n-1}(<\!\! s')$,
$x_i\in k^p\{b_j:j\le s(i)\}$ such that 
$a'\omega_{s'}=v'+d\lln x_2\wedge\cdots\wedge d\lln x_n$.
Thus, 
$a\omega_s=v_1\pm d\lln c\wedge v'\pm d\lln c\wedge d\lln x_2\wedge\cdots\wedge d\lln x_n$.
\qed 
\endpf

\HHH A2.2. Injectivity of the differential symbol 

\phantom{}\smallskip \par 

We can assume that $k$ is a finitely generated field over $\Bbb F_p$.
Then there is a finitely generated algebra over $\Bbb F_p$
with a local ring 
being a discrete valuation ring $\Cal O$ such that
$\Cal O/\Cal M$ is isomorphic to $k$ and the field of fractions $E$ of $\Cal O$ 
is purely transcendental over $\Bbb F_p$.

Using standard results on $K_n(l(t))$ and $\Omega_{l(t)}^n$
one can show that the injectivity of $d_l$ implies the injectivity of $d_{l(t)}$.
Since $d_{\Bbb F_p}$ is injective, so is $d_E$.

Define $k_n(\Cal O)=\kr(k_n(E)\to k_n(k))$.
Then $k_n(\Cal O)$ is generated by symbols
and there is a  homomorphism
$$k_n(\Cal O)\to k_n(k),\quad  \{a_1,\dots,a_n\}\to \{\Overline{a_1},\dots,\Overline{a_n}\},
$$
where $\Overline{a}$ is the residue of $a$.
Let $k_n(\Cal O,\Cal M)$ be its kernel.

Define
$\nu_n({\Cal O})=\kr(\Omega_{\Cal O}^n\to \Omega_{\Cal O}^n/d\Omega_{\Cal O}^{n-1})$,
$\nu_n(\Cal O,\Cal M)=\kr (\nu_n(\Cal O)\to \nu_n(k))$.
There is a homomorphism 
$k_n(\Cal O)\to \nu_n(\Cal O)$ such that
$$\{a_1,\dots, a_n\}\mapsto d\lln a_1\wedge\dots \wedge d\lln a_n.$$

So there is a commutative diagram
$$
\CD
0@>>> k_n(\Cal O,\Cal M) @>>> k_n(\Cal O) @>>> k_n(k) @>>> 0\\
@. @V\varphi VV @VVV @Vd_k VV @.\\
0 @>>> \nu_n(\Cal O,\Cal M) @>>>  \nu_n({\Cal O})@>>> \nu_n(k)    @. .
\endCD
$$
Similarly to A2.1 one can show that $\varphi$ is surjective \cite{BK, 
Prop. 2.4}.
Thus, $d_k$ is injective. \qed

\Bib References

\rf{BK}
S. Bloch and K. Kato, 
{$p$-adic \'etale cohomology},  
Inst. Hautes \'Etudes Sci. Publ. Math.  {63},
(1986), 107--152.

\rf{I} L. Illusie, Complexe de de Rham--Witt et
cohomologie cristalline,
Ann. Sci. \'Ecole Norm. Sup.(4), 12(1979), 501--661.

\rf{K}
K. Kato, 
Galois cohomology of complete discrete valuation fields, 
In {Algebraic $K$-theory}, Lect. Notes in Math. 967, Springer-Verlag, Berlin, 1982, 215--238.

\endBib

\Coordinates

Department of Mathematics \ 
Tokyo Metropolitan University 

Minami-Osawa 1-1, Hachioji, Tokyo 192-03, Japan

E-mail: m-kuri\@comp.metro-u.ac.jp  
\endCoordinates

\smallskip

\Coordinates

Department of Mathematics \  
University of Nottingham

Nottingham NG7 2RD England

E-mail: ibf\@maths.nott.ac.uk
\endCoordinates

\vfill
\pagebreak
\end